\documentclass[twocolumn]{autart}

\usepackage{graphicx}          % Include this line if your
\usepackage{color}
\usepackage{graphicx}
\usepackage{amsfonts,amssymb}
\usepackage{bm}
\usepackage{amsmath}
\usepackage{mathrsfs}
\usepackage{dsfont}
\usepackage{cases}
\usepackage{arydshln}
\usepackage{empheq}
\usepackage{caption}
\usepackage{booktabs} % 用于绘制三线表
\usepackage{tabularx}
\usepackage[table]{xcolor}
\usepackage{dblfloatfix}%表格标题跨栏居中
\usepackage[format=hang]{caption}% 图表标题对齐
\usepackage{float}
\allowdisplaybreaks[4]%%{公式自动分页}
\newtheorem{theorem}{Theorem}[section]

\newtheorem{lemma}{Lemma}[section]

\newtheorem{remark}{Remark}

\newtheorem{assumption}{Assumption}[section]

\begin{document}

\begin{frontmatter}

\title{A Distributed Time-Varying Optimization Approach Based on an Event-Triggered Scheme\thanksref{footnoteinfo}}

\thanks[footnoteinfo]{This paper was supported by the National Natural Science Foundation of China (62176073,12271127), Taishan Scholars of Shandong Province (tsqn202211090), and in part by the Natural Scientific Research Innovation Foundation in Harbin Institute of Technology.}

 \author[1]{Haojin Li}\ead{19S030193@stu.hit.edu.cn},
 \author[3]{Xiaodong Cheng}\ead{xiaodong.cheng@wur.nl},
 \author[3]{Peter van Heijster}\ead{peter.vanheijster@wur.nl},
 \author[1]{Sitian Qin\thanksref{mycorrespondingauthor}}\ead{qinsitian@hitwh.edu.cn}
 % \author[1]{Sitian Qin\corref{mycorrespondingauthor}}
 % \cortext[mycorrespondingauthor]{Corresponding author}
 % \ead{qinsitian@hitwh.edu.cn}

 \address[1]{Department of Mathematics, Harbin Institute of Technology, Weihai, 264209,  China. }  % Please supply
 \address[3]{Mathematical and Statistical Methods Group (Biometris), Wageningen University and Research, Wageningen, 6700 AA, The Netherlands }
\thanks[mycorrespondingauthor]{Corresponding author}


\begin{abstract}
In this paper, we present an event-triggered distributed  optimization approach including a distributed controller to solve a class of distributed time-varying optimization problems (DTOP).  The proposed approach is developed within a distributed neurodynamic (DND) framework that not only optimizes the global objective function in real-time, but also ensures that the states of the agents converge to consensus. This work stands out from existing methods in two key aspects. First, the distributed controller enables the agents to communicate only at designed instants rather than continuously by an event-triggered scheme, which reduces the energy required for agent communication. Second, by incorporating an integral mode technique, the event-triggered distributed controller avoids computing the inverse of the Hessian of each local objective function, thereby reducing computational costs.
Finally, an example of battery charging problem is provided to demonstrate the effectiveness of the proposed event-triggered distributed optimization approach.
\end{abstract}

\end{frontmatter}

\section{Introduction}
Recently, applications such as power systems, transportation, and communication networks, are undergoing a technological transformation since the infrastructural platforms are transforming into complex networked systems with time-varying settings \cite{9133310}. This has spurred research in time-varying optimization which enables decision-making in real-time \cite{DING2025111882,8502778,9133310}.
The practical significance of time-varying optimization compared to static optimization lies in its ability to better describe the real-time properties of some actual engineering tasks.
% As an important optimization problem,  the practical significance of time-varying optimization is in better revealing underlying real-time properties in some actual engineering tasks.
Thus, it has become  prevalent across many engineering applications such as
% optimal control
% for disturbed unmanned helicopter \cite{YAN20241},
  multiquadrotor hose transportation \cite{9670443}, state-of-charge balancing problem \cite{7548310},
%robotic manipulator control \cite{WEN202442},
%formation control \cite{9376701},
 multi-robot navigation problems \cite{9855238} and so on.

It is to be noted that a time-varying cost function means that the optimal solution to such a problem also varies over time and forms a trajectory in its state space. As a result,  traditional algorithms designed for static optimization problems cannot be adopted directly. One popular method for solving time-varying optimization problems is the neurodynamic approach, which is based on recurrent neural networks and capable of performing time-varying optimization in real-time \cite{network1993}.
 Due to its nature of brain-like information processing, neurodynamic approaches can operate in a parallel manner and do not impose strong conditions for global convergence that are difficult to satisfy for most real-world problems \cite{neural1985}. This also implies that it has potential advantages in solving complex optimization problems, especially in time-varying cases.

However, centralized neurodynamic approaches are overly dependent on the information processing power of the central node, limiting their ability to cope with the growing scale of data in industrial systems and scientific problems \cite{8385114}. To address this, distributed neurodynamic (DND)  approaches have attracted a lot of research attention \cite{FIROUZBAHRAMI2022110358,LI2022110259,7862771}. These approaches have proven effective when applied to a distributed time-varying optimization problem (DTOP) \cite{DING2025111882,8100702,7518617,7902101,7862771}. For example, a DND approach was proposed in \cite{DING2025111882} for the networked Lagrangian agents, and zero optimum-tracking error was realized. In \cite{7902101}, a  discrete DND approach was proposed based on the prediction-correction methods. In \cite{7862771}, a gradient-based searching method was combined with a DND algorithm to solve a DTOP with quadratic cost functions. In \cite{7518617}, two DND algorithms were designed to solve DTOP with a general local objective function. In \cite{8100702}, a DND algorithm based on a fixed-time consensus approach was proposed.  However, in above studies, all the  agents have to communicate with each other continuously for real-time updates of their states, which inevitably leads to the over-consumption of communication resources.

To address this issue, we introduce an event-triggering scheme that initiates communication only when necessary.
Event-triggering schemes are commonly considered in network control to enhance resource efficiency \cite{CAI2024111508,HUANG2024111877,9640505,WU202496}, shifting the conventional continuous control to a conscious intermittent control \cite{WU202496}.  For instance, an event-triggered approach was proposed in \cite{9640505}, where each agent broadcasts its state only if the error between the current state and the last broadcast state is bigger than the given upper bound. In \cite{WANG201834}, two types of event-triggered schemes were investigated for a double-integrator multi-agent system, resulting in reduced communication consumption. However, most of the work combines event-triggered schemes with static optimization algorithms. To the best of our knowledge, there is no research that has integrated event-triggered schemes into DTOPs.  Hence, how to design distributed algorithms with intermittent communication to solve DTOPs is still an open question.

Besides, for algorithm design, the aforementioned DND algorithms, and those in \cite{10339857,10402058,9855238}, typically involve computing the inverse of the Hessian of the local objective function. This requirement imposes significant computation demands as it is inherently computationally expensive to compute the inverse of a matrix, particularly for high-dimensional systems.
%This often demands a high computational load, particularly for high-dimensional systems.
Therefore, it is both theoretically and practically important to explore methods to avoid computing the inverse of matrices in the optimization algorithms. To the best of our knowledge, only a few solutions have been proposed to consider this issue in the centralized framework (see e.g., \cite{CHEN2024106462,10003122}). Addressing this issue becomes even more challenging within the event-triggered and distributed framework of this paper.

In this paper, we address the aforementioned challenges in the DTOPs. Specifically, an event-triggered scheme is adopted in designing a distributed controller, aiming to solve a DTOP without requiring continuous-time communication among the distributed agents. Our main contributions are highlighted as threefold.

\begin{enumerate}
\item For the first time, an event-triggered scheme is incorporated with the distributed neurodynamic optimization approach design,  leading to a communication-efficient distributed controller. Compared to the algorithms for continuous-time communication in e.g., \cite{10339857,10402058,8100702,9855238}, the distributed neurodynamic optimization approach in this paper has the potential to save communication resources.

\item A computationally  efficient distributed neurodynamic optimization approach is proposed. Compared to e.g. \cite{10339857,10402058,8100702,9855238}, where the inverse of the Hessian and partial time derivative of the gradient of the local objective function are required,  the proposed approach only requires the information of the Hessian matrix and the gradient to be transferred among the agents, which can significantly reduce the computational cost.
%This information can potentially provide control engineers with a priori information about the system’s behavior.
\item Compared to the works in e.g.,\cite{10339857,doi:10.1080/00207721.2020.1801885} based on estimators, our approach requires fewer variables and agent interactions. We also remove the need for the boundedness of the second-order time derivative of the gradient function and the time derivative of the Hessian matrices as required in \cite{doi:10.1080/00207721.2020.1801885}.
This  relaxes the implementation conditions of the distributed neurodynamic optimization approach.
\end{enumerate}

This paper is structured as follows: In Section \ref{sec2}, we give the problem formulation and state some preliminary knowledge and assumptions on the problem. In Section \ref{sec3}, the  event-triggered DND approach is proposed to solve a class of DTOPs, and the consensus and convergence properties are analyzed. In Section \ref{sec4}, an example related to state-of-charge balancing problem is given to show the effectiveness of the approach. The conclusions are given in Section \ref{sec5}.
%\vspace{-3mm}
\section{Preliminaries}\label{sec2}
\subsection{Notations}
Let $\mathbb{R},\mathbb{R}^n, \mathbb{R}^{m\times n}, \mathbb{R}_+$ be the sets of real numbers, $n$-dimension real vectors, $m\times n$ real matrices and nonnegative real numbers, respectively. For $A\in \mathbb{R}^{m\times n}$, denote $A^\mathrm{T}$ as the transpose of $A$. Denote $\mathbf{1}_N$ (or  $\mathbf{0}_N$) as the $N$-dimensional vector with all elements of 1 (or 0). $\|\cdot\|_1$ and $\|\cdot\|$ represent the 1-norm and the Euclidean norm. Furthermore,  $\nabla_xf(x,t)$, $\nabla^2_x f(x,t)$ are defined as the partial derivative and the Hessian matrix of $f(x,t)$ with respect to $x$, while $\frac{\partial}{\partial t}\nabla_xf(x,t)$ is the partial derivative of $\nabla_xf(x,t)$ with respect to $t\in \mathbb{R}_+$. $I_n$ represents the $n\times n$ identity matrix. For a vector $v = (v_1, v_2, ... , v_n)^\mathrm{T}\in \mathbb{R}^n$,  $\mathrm{sign}(v)=(\mathrm{sign}(v_1),\mathrm{sign}(v_2),\dots,\mathrm{sign}(v_n))^\mathrm{T}$ with $\mathrm{sign}(v_i)$ the signum function defined as
\begin{equation*}
\mathrm{sign}(v_i)=\left\{\begin{split}
&1, \qquad &&v_i>0\\
&0, \qquad &&v_i=0\\
&-1, \qquad &&v_i<0.
\end{split}\right.
\end{equation*}

\subsection{Graph Theory}
The interconnection topology of a multi-agent system can be represented by a bidirectional graph $\mathcal{G}=(\mathcal{V},\mathcal{E})$, where the node set is $\mathcal{V}=\{1,2,\dots,N\}$ and the edge set is $\mathcal{E}\subset\mathcal{V}\times\mathcal{V}$. The adjacency matrix of a graph $\mathcal{G}$, denoted by $\mathcal{A}=[a_{ij}]\in \mathbb{R}^{N\times N}$, is defined such that $a_{ij}=1$ if $(j,i)\in \mathcal{E}$  and $a_{ij}=0$ otherwise. The graph is bidirectional if and only if $a_{ij}=a_{ji}$ for all $i,j \in \mathcal{V}$. The set $\mathcal{N}_i=\{j\in\mathcal{V}: (j,i)\in\mathcal{E}\}$ collects all the neighbors of the $i$th node.
In a graph $\mathcal{G}$, a path  between node $i_1$ and node $i_k$ is a sequence of edges of the form $(i_1,i_2),(i_2,i_3),\dots, (i_{k-1},i_k)$ $i_k\in \mathcal{V}$.
If there exists a path between every pair of nodes, then the graph is considered connected.

% Let $\mathcal{L}=[l_{ij}]\in\mathbb{R}^{N\times N}$ be the Laplacian matrix associated with $\mathcal{A}$, in which $l_{ii}=\sum_{j=1,j\neq i}^Na_{ij}$ and $l_{ij}=-a_{ij}$ for  $i\neq j$.
% For a graph $\mathcal{G}$ which is bidirectional and connected, the Laplacian matrix $\mathcal{L}$ is symmetric positive semidefinite satisfying $\mathcal{L}\mathbf{1}_N=\mathbf{1}_N^\mathrm{T}\mathcal{L}=\mathbf{0}_N$.

% \begin{lemma}\label{le_en}(Lemma 1 in \cite{9359452})
% For a graph $\mathcal{G}$ which is bidirectional and connected, the Laplacian matrix $\mathcal{L}$ is symmetric positive semidefinite satisfying $\mathcal{L}\mathbf{1}_N=\mathbf{1}_N^\mathrm{T}\mathcal{L}=\mathbf{0}_N$. For a vector $x\in \mathbb{R}^N$,  $x^\mathrm{T}\mathcal{L}x=(1/2)\sum_{i=1}^N\sum_{j=1}^Na_{ij}(x_i-x_j)^2$. For its  eigenvalues, there is the order of $0=\lambda_1(\mathcal{L})<\lambda_2(\mathcal{L})\leq\dots\leq\lambda_N(\mathcal{L})$. Here, $\lambda_2(\mathcal{L})$ represents the minimum nonzero eigenvalue of the Laplacian matrix $\mathcal{L}$ and $\lambda_2(\mathcal{L})=\min_{x\neq\mathbf{0}_N,\mathbf{1}^\mathrm{T}_Nx=0}(x^\mathrm{T}\mathcal{L}x/x^\mathrm{T}x)$.
% \end{lemma}

\subsection{Problem Description}
Consider the following time-varying optimization problem:
\begin{equation}\begin{split}\label{eq2}
\min_{x(t)} ~~&F(x(t),t),
\end{split}\end{equation}
 where $F(x(t),t): \mathbb{R}^n\times \mathbb{R}_+\rightarrow \mathbb{R}$.
 Solving  \eqref{eq2} using centralized optimization algorithms is often computationally expensive, particularly when the time-varying objective function $F(x(t),t)$ is complex. To reduce the computational load, we can implement distributed computation with a multi-agent system, provided that $F(x(t),t)$ can be rewritten as the sum of individual time-varying functions as $$F(x(t),t)=\sum_{i=1}^Nf_i(x(t),t).$$
 In this setup,  each agent $i$ possesses a local state $x_i(t)$ and is responsible for optimizing a local objective function $f_i(x_i(t),t)$, which is only known to the agent itself.
 %Notably, the centralized optimization algorithms for solving \eqref{eq2} frequently necessitates a central agent with powerful computational capabilities and all the other agents share complete information with each other. In order to ease the burden on the central agent and private protection,  we assign each agent a local objective function $f_i(x_i(t),t)$  which is only known to the $i$th agent.
 The $N$ agents can communicate with each other, and the communication network topology is represented by a bidirectional connected graph $\mathcal{G}$.

Then, the state of each agent is designed to
evolve based on the following dynamics:
\begin{equation}\label{agent_dynamic}
\dot{x}_i(t)=u_i(t),
\end{equation}
where $u_i(t)\in\mathbb{R}^n$ is the control input of the $i$th agent to be specified.
With \eqref{agent_dynamic},  the problem \eqref{eq2} can be transformed into a distributed time-varying optimization problem (DTOP) where all the agents achieve \textit{consensus} while optimizing the team objective function  $\sum_{i=1}^Nf_i(x_i(t),t)$. This leads to the following formulation:
% Notably, the design of the approach for such a problem \eqref{eq2} is significantly increased since the optimization of centralized issues frequently necessitates a central agent with powerful computational capabilities. In order to ease the burden on the central agent and private protection,  we assign each agent an objective function $f_i(x,t)$ (called local objective function) when $F(x,t)$ can be rewritten as a sum of local objective functions ($F(x,t)=\sum_{i=1}^Nf_i(x,t)$) and the TVOP \eqref{eq2} can be transformed  as
\begin{equation} \label{DTOP}
\begin{aligned}
    \min_{x_i(t)} ~~&\sum_{i=1}^Nf_i(x_i(t),t)\\
\mathrm{s.t.}~~~&x_i(t)=x_j(t),~~ \forall i,j\in \mathcal{V},
\end{aligned}
\end{equation}
where $x_i(t)\in\mathbb{R}^n$ is the state of the $i$th agent.  The goal becomes to design $u_i(t)$ in \eqref{agent_dynamic} for each agent $i$  by only exploiting its local information from the interaction with its neighbors such that all the agents can cooperatively track an optimal trajectory $x^*(t)$ of the DTOP \eqref{DTOP} (The optimal trajectory is unique when Assumption \ref{strongly_convex} is made). For notational simplicity, throughout the paper, we remove the time index $t$ from the variable $x_i(t)$ and $u_i(t)$ and only keep it when necessary. Additionally, the following assumptions are made for the theoretical analysis in our result.
%   \begin{assumption}\label{a1}
% The fixed graph $\mathcal{G}$ is bidirectional and connected.
% \end{assumption}
\begin{assumption}\label{strongly_convex}
    For each $i$, the local objective function $f_i(x,t)$ is twice continuously differentiable and strongly convex with respect to $x$ and continuously differentiable with respect to $t$.
\end{assumption}
\begin{assumption}\label{as3}
There exists a $\underline{\lambda}_i>0$, such that $\underline{\lambda}_i\leq\min\{\lambda(\nabla^2_{x_i}f_i(x,t))|x\in \mathbb{R}^n, t\geq 0\}$.
% $\underline{\lambda}_{\min}\triangleq\inf_{t\geq 0}\{\lambda_{\min}(t)\}>0$, where $\lambda_{\min}(t)$ is the minimum eigenvalue of $\nabla_{x_i}^2f_i(x_i,t)$.
\end{assumption}

\begin{assumption}\label{as2}
There exists a $\overline{f}>0$ such that $$\left\|\frac{\partial}{\partial t}\nabla_{x}f_i(x,t)\right\|\leq\overline{f}$$
holds for all $i\in\mathcal{V}$ and $x\in \mathbb{R}^n, t\geq0$.
\end{assumption}

\begin{remark}
Assumptions \ref{as3} and \ref{as2} are necessary for the theoretical analysis on the convergence of the proposed method, as detailed in the proof of Lemma~\ref{le_op_fix}.  Note that Assumption \ref{as3} ensures that $\nabla_{x_i}^2f_i(x_i,t)$  is positive definite for all $t\geq 0$. Similar conditions are  required in previous works, e.g., \cite{10402058,9855238}.
 Furthermore,
   Assumption~\ref{as2} is also adopted in distributed time-varying optimization, e.g., \cite{10402058,7518617}.
 An example of objective functions satisfying Assumptions \ref{as3} and \ref{as2} is: $f_i(x,t)=(a_ix+b_i(t))^2$, which is commonly used in energy minimization with $\| \dot{b}_i(t) \|$ upper  bounded for all $i\in\mathcal{V}$.
 % It holds for a vital kind of functions. For instance,  the objective function commonly in energy minimization $f_i(x,t)=(a_ix+b_i(t))^2$. It is clear that Assumption \ref{as2} holds as long as $\dot{b}_i(t)$  is bounded for all $i\in\mathcal{V}$.
\end{remark}

In this paper, we aim to overcome the limitation in the current literature on DTOPs where each agent in the network has to continuously communicate with the others and update its state. To minimize communication resource usage, we propose to use an event-triggering scheme,
 % which has been addressed in the other problems, see e.g., \cite{9669076,9640505,9657491,9709887,9448297,8949816,7934362}.
where the distributed agents only communicate when an
event-triggering condition is satisfied. Then, the research questions to be addressed in this paper are: How to design an event-triggered controller in \eqref{agent_dynamic} to solve the DTOP in \eqref{DTOP}? Can we theoretically guarantee the convergence of each agent to the optimal solution?

\section{Main Results}\label{sec3}
In this section, we present a novel method to solve the DTOP problem in \eqref{DTOP}, which integrates the event-triggered scheme into the distributed neurodynamic (DND) approach, and we provide a theoretical analysis of the convergence of the proposed method.

By adapting the DND controller structure in \cite{10402058}, we propose the following distributed controller with event triggering mechanism for the multi-agent system \eqref{agent_dynamic}:
\begin{subnumcases}{\label{controller}}
u_i=-k_1s_i-k_2\sum_{j\in \mathcal{N}_i}\mathrm{sign}(\widetilde{x}_i-\widetilde{x}_j)\label{controller1}\\
~~~~~~~\notag-k_3\mathrm{sign}(\nabla^2_{x_i}f(x_i,t)s_i),\\
s_i=
\nabla_{x_i}f_i(x_i,t)\label{controller2}\\
~~~~~~~\notag+k_2\int_{0}^t\nabla^2_{x_i}f_i(x_i,\tau)\sum_{j\in\mathcal{N}_i}\mathrm{sign}(\widetilde{x}_i-\widetilde{x}_j)\mathrm{d}\tau,
\end{subnumcases}
for $t\in [t_k^i,t_{k+1}^i)$, where $t_k^i$ is the $k$th event-triggering instant of agent $i$, and  $\widetilde{x}_i(t)=x_i(t_k^i)$, for all
$t\in[t_k^i,t_{k+1}^i)$, is defined as the event-triggering state of $x_i(t)$. Similarly, $\widetilde{x}_j(t)=x_j(t_{k'}^j)$ for all $t\in[t_{k'}^j,t_{k'+1}^j)$ is the event-triggered state of $x_j(t)$.
% with $t_{k'}^j$ denoting the latest triggering time for agent $j$ before the event of agent $i$ is triggered.
The tuning parameters $k_1,k_2,k_3$ are all positive, and $s_i(t)$ is an auxiliary variable designed to track the gradient sum $\sum_{i=1}^N\nabla_{x_i}f_i(x_i,t)$. Specifically, as $s_i \rightarrow 0$ for each $i$  the gradient sum also converges to zero.
%and it helps $\sum_{i=1}^N\nabla_{x_i}f_i(x_i,t)$  approach to zero.

\begin{figure}[!h]
   \begin{center}
       \includegraphics[scale=0.3]{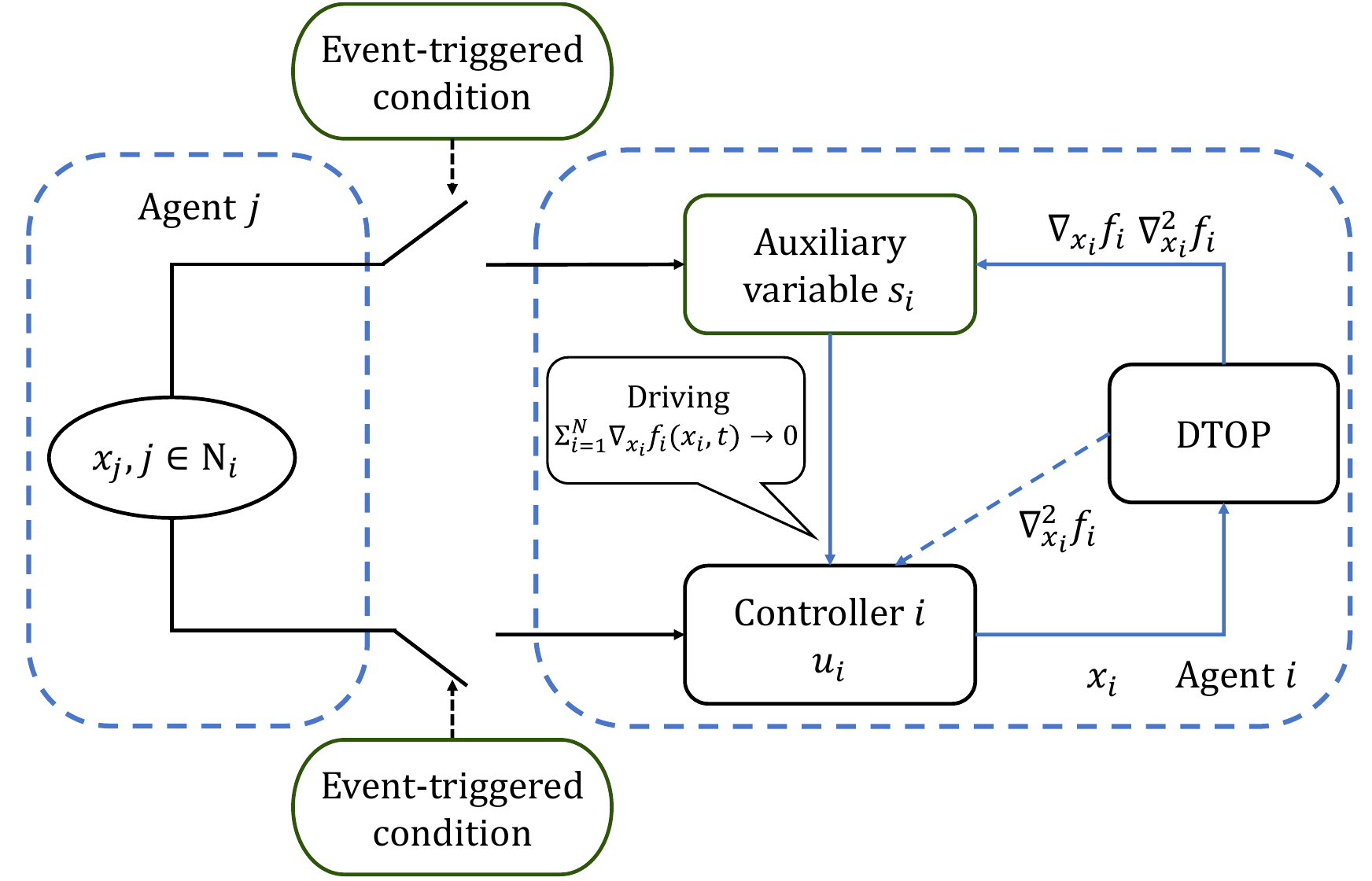}
        \caption{Block diagram of the proposed event-triggered DND approach. }
        \label{fig_diagram}
        \end{center}
\end{figure}

 Inspired by \cite{9359452}, for each agent $i$, an event-triggered scheme is designed as:
\begin{equation}\label{trigger_scheme}
t_{k+1}^i=\inf\limits_{t>t_k^i}\{t: G_{i}(t)>0\},
\end{equation}
with $G_{i}(t)$ the triggering function defined as:
\begin{equation}\label{trigger_function}
\begin{split}
% G_{i}(t)&=||\varepsilon_{i}(t)||-\frac{1}{2N}\sum_{j\in\mathcal{N}_i}||\widetilde{x}_i-\widetilde{x}_j||-m_i\mathrm{e}^{-a_it},\\
G_{i}(t)&=\|\varepsilon_{i}(t)\|-m_i\mathrm{e}^{-a_it},\\
\end{split}\end{equation}
where  $\varepsilon_{i}(t)$  represents the  measurement error of agent $i$:
\begin{align}\label{measurement_error}
\varepsilon_{i}(t)=\widetilde{x}_i(t)-x_i(t).
\end{align}
The parameters $m_i,a_i$ are all positive and  the term $m_i\mathrm{e}^{-a_it}$ in the triggering function \eqref{trigger_function} can effectively prevent the Zeno behavior \cite{9359452}.
%in which the Zeno behavior refers to
That is, the
triggering of an infinite number
of events within a limited period of time. The structural block diagram of the proposed event-triggered DND approach is displayed in Fig. \ref{fig_diagram}.

In the event-triggered scheme \eqref{trigger_scheme}, if the triggering condition $G_{i}>0$ is fulfilled, then $x_i(t_k^i)$  updates its current state while broadcasting its current state information to its neighbors, and $\varepsilon_{i}$  will thus be reset to 0. With the triggering function \eqref{trigger_function}, the agents communicate only at the designed instances, by which the communication frequency can be significantly reduced.
% Thus, compared with the existing distributed time-varying optimization algorithms in e.g., \cite{10339857,10402058,8100702,9855238}, less communication resource is required.
To the best of our knowledge, there are no studies on integrating event-triggered schemes into DTOPs in the existing literature.

\begin{remark}
Compared to
    the distributed time-varying optimization algorithms in \cite{10339857,10402058,8100702,9855238},
    % require the information of  $\nabla_{x_i}f_i(x_i,t)$, $(\nabla_{x_i}^2f_i(x_i,t))^{-1}$ and $\frac{\partial}{\partial t}\nabla_{x_i}f_i(x_i,t)$, in distributed controller $u_i$,
    we remove the requirement for information on $(\nabla_{x_i}^2f_i(x_i,t))^{-1}$ and $\frac{\partial}{\partial t}\nabla_{x_i}f_i(x_i,t)$ in  the distributed controller \eqref{controller}. Hence, the computation load is reduced. In particular, the computation of the inverse of a matrix requires a lot of computational resources, especially in some large-scale matrices.  Therefore, the distributed controller \eqref{controller} is more suitable for solving complex problems than the algorithms in \cite{10339857,10402058,8100702,9855238}. Besides, compared with estimator-based average tracking methods \cite{10339857,doi:10.1080/00207721.2020.1801885}, the proposed approach has fewer variables and
communication load among the agents, contributing to a lower computational cost.
\end{remark}

% \begin{remark}\label{rem_scheme}

% \end{remark}

The following theorem guarantees the consensus and convergence of the distributed controller \eqref{controller}.
\begin{theorem}\label{th_fixtop}
Consider the distributed controller \eqref{controller} with the event-triggered scheme \eqref{trigger_scheme} for each agent $i$. Under Assumptions \ref{as3}-\ref{as2}, if
\begin{subequations}
\label{k_3}
    \begin{align}
k_3&>\overline{f}/\underline{\lambda}_i \label{K3},
\\
\nabla^2_{x_i}f_i(x_i,t)&=\nabla^2_{x_j}f_j(x_j,t),~~ \forall i,j\in \mathcal{V}, t\geq 0, \label{hessian}
 \end{align}
\end{subequations}
 then the state of each agent $x_i$  reaches consensus and  cooperatively track the optimal trajectory $x^*(t)$ of the DTOP \eqref{DTOP}. Furthermore, the Zeno behavior is excluded.
\end{theorem}

In order to prove Theorem \ref{th_fixtop}, we will present three lemmas in the following. The result of Theorem \ref{th_fixtop} then follows immediately from the three lemmas. In Lemma~\ref{le_op_fix}, we prove that the
distributed controller \eqref{controller2} ensures $\sum_{i=1}^N\nabla_{x_i}f_i(x_i,t)\rightarrow\mathrm{0}_n$ in a finite time, and Lemma \ref{le_consensus_fix} then shows that the states of all the agents eventually reach consensus on the optimal trajectory. These two lemmas guarantee that the DTOP in \eqref{DTOP} can be solved by using the event-triggered controller \eqref{controller}-\eqref{measurement_error}. Finally, Lemma  \ref{le_zeno} shows that the Zeno behavior can be avoided by our method.

\begin{lemma}\label{le_op_fix}
Assume the conditions in Theorem~\ref{th_fixtop} hold,
%Under  Assumption \ref{a1}-\ref{as2}, if condition \eqref{k_3} holds
for the system \eqref{agent_dynamic} with the distributed controller \eqref{controller}, there is a constant $T_{\max}>0$ such that $\sum_{i=1}^N\nabla_{x_i}f_i(x_i,t)=0$, for all $t\geq T_{\max}$.
\end{lemma}
\begin{pf}
Consider the following Lyapunov function for the agent $i$:
$$V_{1i}=\frac{1}{2}s_i^\mathrm{T}s_i,$$
which has the time derivative  along \eqref{controller2} as
\begin{align*}
\dot{V}_{1i}
% =&s_i^\mathrm{T}\dot{s}_i
% \\
=&s_i^\mathrm{T}\big[\nabla^2_{x_i}f_i(x_i,t)\dot{x}_i+\frac{\partial}{\partial t}\nabla_{x_i}f_i(x_i,t)\\
&+k_2\nabla^2_{x_i}f_i(x_i,t)\sum_{j\in \mathcal{N}_i}\mathrm{sign}(\widetilde{x}_i-\widetilde{x}_j)\big]
\\
% =&s_i^\mathrm{T}\big[\nabla^2_{x_i}f_i(x_i,t)\big(-k_1s_i-k_2\sum_{j\in\mathcal{N}_i}\mathrm{sign}(\widetilde{x}_i-\widetilde{x}_j)\\
% &-k_3\mathrm{sign}(\nabla^2_{x_i}f_i(x_i,t)s_i)\big)\\
% &+\frac{\partial}{\partial t}\nabla_{x_i}f_i(x_i,t)\\
% &+k_2\nabla^2_{x_i}f_i(x_i,t)\sum_{j\in\mathcal{N}_i}\mathrm{sign}(\widetilde{x}_i-\widetilde{x}_j)\big]\\
=&-k_1s_i^\mathrm{T}\nabla^2_{x_i}f_i(x_i,t)s_i
+s_i^\mathrm{T}\frac{\partial}{\partial t}\nabla_{x_i}f_i(x_i,t)\\
%&-k_2s_i^\mathrm{T}\nabla^2_{x_i}f_i(x_i,t)\sum_{j\in\mathcal{N}_i}\mathrm{sign}(\widetilde{x}_i-\widetilde{x}_j)\\
&-k_3s_i^\mathrm{T}\nabla^2_{x_i}f_i(x_i,t)\mathrm{sign}(\nabla^2_{x_i}f_i(x_i,t)s_i),
%&+s_i^\mathrm{T}\frac{\partial}{\partial t}\nabla_{x_i}f_i(x_i,t)\\
%&+k_2s_i^\mathrm{T}\nabla^2_{x_i}f_i(x_i,t)\sum_{j\in\mathcal{N}_i}\mathrm{sign}(\widetilde{x}_i-\widetilde{x}_j).
\end{align*}
where the terms containing $k_2$ are canceled out. Note that
$$ s_i^\mathrm{T}\nabla^2_{x_i}f_i(x_i,t)\mathrm{sign}(\nabla^2_{x_i}f_i(x_i,t)s_i) = \|\nabla^2_{x_i}f_i(x_i,t)s_i\|_1. $$
Due to the condition \eqref{hessian}, denote $\underline{\lambda}=\underline{\lambda}_i=\underline{\lambda}_j$.
Then, with Assumption \ref{as3} and Assumption \ref{as2}, we have
% \begin{align*}
% \dot{V}_{1i}&\leq-k_1\underline{\lambda}^i_{min}s_i^\mathrm{T}s_i+s_i^\mathrm{T}\frac{\partial}{\partial t}\nabla_{x_i}f_i(x_i,t)\\
% &\quad -k_3\|\nabla^2_{x_i}f_i(x_i,t)s_i\|_1\\
% &\leq -2k_1\underline{\lambda}^i_{min}V_{1i}+\overline{f}\|s_i\|-k_3\|\nabla^2_{x_i}f_i(x_i,t)s_i\|\\
% &\leq -2k_1\underline{\lambda}^i_{min}V_{1i}+\overline{f}\|s_i\|-k_3\underline{\lambda}^i_{min}\|s_i\|\\
% &=-2k_1\underline{\lambda}^i_{min}V_{1i}-2^\frac{1}{2}(k_3\underline{\lambda}^i_{min}-\overline{f})V_{1i}^\frac{1}{2}\\
% &\leq-2^\frac{1}{2}(k_3\underline{\lambda}^i_{min}-\overline{f})V_{1i}^\frac{1}{2}.
% \end{align*}
\begin{align*}
\dot{V}_{1i}&\leq  -k_1\underline{\lambda}s_i^\mathrm{T}s_i+\overline{f}\|s_i\|-k_3\|\nabla^2_{x_i}f_i(x_i,t)s_i\|\\
&\leq -2k_1\underline{\lambda}V_{1i}+\overline{f}\|s_i\|-k_3\underline{\lambda}\|s_i\|\\
&=-2k_1\underline{\lambda}V_{1i}-2^\frac{1}{2}(k_3\underline{\lambda}-\overline{f})V_{1i}^\frac{1}{2}\\
&\leq-2^\frac{1}{2}(k_3\underline{\lambda}-\overline{f})V_{1i}^\frac{1}{2}.
\end{align*}

%Barbalat Lemma \cite{slotine1991applied} and
Because of $k_3>\overline{f}/\underline{\lambda}$ and  Lemma 2.7 in \cite{LI2024127589}, there exists a  constant $T^i_{\max}>0$
upper bounded by $$T^i_{\max}\leq\frac{\sqrt{2}\sqrt{V_{1i}(0)}}{k_3\underline{\lambda}-\overline{f}}$$
such that $V_{1i}=0$ when $t>T^i_{\max}$.
Denote $T_{\max}=\max_i\{T^i_{\max}\}$. Hence, $s_i=0$ for all $i\in\mathcal{V}$ when $t>T_{\max}$. Considering the $\mathcal{G}$ is bidirectional and connected, $\nabla^2_{x_i}f_i(x_i,t)=\nabla^2_{x_j}f_j(x_j,t)$ for all $i,j\in \mathcal{V}, t\geq 0$, by summing \eqref{controller2}, it can be concluded that $\sum_{i=1}^Ns_i=\sum_{i=1}^N\nabla_{x_i}f_i(x_i,t)=0$ when $t>T_{\max}$.
$\hfill{} \Box$
\end{pf}

\begin{lemma}\label{le_consensus_fix}
Assume the conditions in Theorem \ref{th_fixtop} hold,
for system \eqref{agent_dynamic} with the distributed controller \eqref{controller}, all the state $x_i$ will reach consensus asymptotically on the optimal trajectory $x^*(t)$ of the DTOP \eqref{DTOP}.
\end{lemma}
\begin{pf}
Consider $T_{\max}$ in Lemma~\ref{le_op_fix}, which leads to $\sum_{i=1}^N\nabla_{x_i}f_i(x_i,t)=0$ when  $t>T_{\max}$, and then the distributed controller \eqref{controller1} is simplified as
\begin{align}\label{recontroll}
u_i=-k_2\sum_{j\in\mathcal{N}_i}\mathrm{sign}(\widetilde{x}_i-\widetilde{x}_j).
\end{align}
Define the  consensus error for the agent $i$ as
 \begin{align}\label{consensus_error}
\widehat{x}_i : =x_i-\frac{1}{N}\sum_{j=1}^Nx_j,
 \end{align}
 % It follows from \eqref{recontroll} that we can get
 % \begin{equation}\label{dynamic_error}
 %     \dot{\widehat{x}}_i=-k_2\widehat{x}_i-k_4\sum_{j\in\mathcal{N}_i}\mathrm{sign}(\widetilde{x}_i-\widetilde{x}_j)
 % \end{equation}
 and construct a Lyapunov function as
\begin{align*}
V_{2}=\frac{1}{2}\sum_{i=1}^N\widehat{x}_i^\mathrm{T}\widehat{x}_i+k_2N(1+N^\frac{1}{2})\sum_{i=1}^N\frac{m_i}{a_i}\mathrm{e}^{-a_it}.
\end{align*}

Due to $\sum_{i=1}^N\widehat{x}_i=0$,  the time derivative of $V_{2i}$ can be obtained as
\begin{align*}
\dot{V}_{2}=&-k_2\sum_{i=1}^N\sum_{j\in\mathcal{N}_i}\widehat{x}_i^\mathrm{T}\mathrm{sign}(\widetilde{x}_i-\widetilde{x}_j)\\
&-k_2N(1+N^\frac{1}{2})\sum_{i=1}^Nm_i\mathrm{e}^{-a_it}.
\end{align*}
Since the underlying graph $\mathcal{G}$ is connected and bidirectional, we have
\begin{align*}
    \sum_{i=1}^N\sum_{j\in \mathcal{N}_i}\widehat{x}_i^\mathrm{T}\mathrm{sign}(\widetilde{x}_i-\widetilde{x}_j)
&=\sum_{i=1}^N\sum_{j=1}^Na_{ij}\widehat{x}_i^\mathrm{T}\mathrm{sign}(\widetilde{x}_i-\widetilde{x}_j)
\\
& =\sum_{j=1}^N\sum_{i=1}^Na_{ji}\widehat{x}_j^\mathrm{T} \mathrm{sign}(\widetilde{x}_j-\widetilde{x}_i).
\end{align*}
Hence,
% \begin{align*}
%     \sum_{i=1}^N\sum_{j\in \mathcal{N}_i}\widehat{x}_i^\mathrm{T}\mathrm{sign}(\widetilde{x}_i-\widetilde{x}_j)= \frac{1}{2}\sum_{i=1}^N\sum_{j\in \mathcal{N}_i}(\widehat{x}_i-\widehat{x}_j)^\mathrm{T}\mathrm{sign}(\widetilde{x}_i-\widetilde{x}_j)
% \end{align*}.
% by exchange the subscripts $i,j$, $\sum_{i=1}^N\sum_{j\in \mathcal{N}_i}\widehat{x}_i^\mathrm{T}\mathrm{sign}(\widetilde{x}_i-\widetilde{x}_j)
% =\sum_{i=1}^N\sum_{j=1}^Na_{ij}\widehat{x}_i^\mathrm{T}\mathrm{sign}(\widetilde{x}_i-\widetilde{x}_j)
% =\sum_{j=1}^N\sum_{i=1}^Na_{ji}\widehat{x}_j^\mathrm{T}\\ \mathrm{sign}(\widetilde{x}_j-\widetilde{x}_i).
% $ Because the graph is bidirectional, $a_{ij}=a_{ji}$, hence, $\sum_{i=1}^N\sum_{j\in \mathcal{N}_i}\widehat{x}_i^\mathrm{T}\mathrm{sign}(\widetilde{x}_i-\widetilde{x}_j)=-\sum_{i=1}^N\sum_{j=1}^Na_{ij}\widehat{x}_j^\mathrm{T}\mathrm{sign}(\widetilde{x}_i-\widetilde{x}_j)=\frac{1}{2}\sum_{i=1}^N\sum_{j\in \mathcal{N}_i}(\widehat{x}_i-\widehat{x}_j)^\mathrm{T}\mathrm{sign}(\widetilde{x}_i-\widetilde{x}_j)$.
$\dot{V}_{2}$ can be rewritten as
\begin{align*}
\dot{V}_{2}=&-\frac{k_2}{2}\sum_{i=1}^N\sum_{j\in\mathcal{N}_i}\big(\widehat{x}_i-\widehat{x}_j\big)^\mathrm{T}\mathrm{sign}(\widetilde{x}_i-\widetilde{x}_j)\\
&-k_2N(1+N^\frac{1}{2})\sum_{i=1}^Nm_i\mathrm{e}^{-a_it}\\
=&-\frac{k_2}{2}\sum_{i=1}^N\sum_{j\in\mathcal{N}_i}\big(x_i-x_j\big)^\mathrm{T}\mathrm{sign}(\widetilde{x}_i-\widetilde{x}_j)\\
&-k_2N(1+N^\frac{1}{2})\sum_{i=1}^Nm_i\mathrm{e}^{-a_it},
% =&-\frac{k_2}{2}\sum_{i=1}^N\sum_{j\in\mathcal{N}_i}\big((x_i-x_j)+(\varepsilon_i-\varepsilon_j)\\
% &-(\varepsilon_i-\varepsilon_j)\big)^\mathrm{T}\mathrm{sign}(\widetilde{x}_i-\widetilde{x}_j)\\
% &-k_2N(1+N^\frac{1}{2})\sum_{i=1}^Nm_i\mathrm{e}^{-a_it}.
\end{align*}
where $\widehat{x}_i-\widehat{x}_j=x_i-x_j$ is applied to obtain the second equation.

It follows from \eqref{measurement_error} that $x_i-x_j=(\widetilde{x}_i-\widetilde{x}_j) - (\varepsilon_i-\varepsilon_j)$, which leads to
\begin{align*}
\dot{V}_{2}=&
-\frac{k_2}{2}\sum_{i=1}^N\sum_{j\in\mathcal{N}_i}\big(\widetilde{x}_i-\widetilde{x}_j\big)^\mathrm{T}\mathrm{sign}(\widetilde{x}_i-\widetilde{x}_j)\\
&+\frac{k_2}{2}\sum_{i=1}^N\sum_{j\in\mathcal{N}_i}\big(\varepsilon_i-\varepsilon_j\big)^\mathrm{T}\mathrm{sign}(\widetilde{x}_i-\widetilde{x}_j)\\
&-k_2N(1+N^\frac{1}{2})\sum_{i=1}^Nm_i\mathrm{e}^{-a_it}\\
\leq&-\frac{k_2}{2}\sum_{i=1}^N\sum_{j\in\mathcal{N}_i}\|\widetilde{x}_i-\widetilde{x}_j\|
+\frac{k_2}{2}N^\frac{1}{2}\sum_{i=1}^N\sum_{j\in\mathcal{N}_i}\|\varepsilon_i-\varepsilon_j\|\\
&-k_2N(1+N^\frac{1}{2})\sum_{i=1}^Nm_i\mathrm{e}^{-a_it}\\
% \end{align*}
% Since $||\widetilde{x}_i-\widetilde{x}_j||=||(x_i-x_j)+(\varepsilon_i-\varepsilon_j)||\geq \\||x_i-x_j||-||\varepsilon_i-\varepsilon_j||$, it can be concluded that
% \begin{align*}
% \dot{V}_{2i}
% \leq&-k_2\sum_{i=1}^N\widehat{x}_i^\mathrm{T}\widehat{x}_i
% -\frac{k_4}{2}\sum_{i=1}^N\sum_{j\in\mathcal{N}_i}||x_i-x_j||_2\\
% &+k_4\sum_{i=1}^N\sum_{j\in\mathcal{N}_i}||\varepsilon_i-\varepsilon_j||_2
% -2k_4N\sum_{i=1}^N m_i\mathrm{e}^{-a_it}\\
\leq&-\frac{k_2}{2}\sum_{i=1}^N\sum_{j\in\mathcal{N}_i}\|x_i-x_j\|\\
&+\frac{k_2}{2}(1+N^\frac{1}{2})\sum_{i=1}^N\sum_{j=1}^N\|\varepsilon_i-\varepsilon_j\|\\
&-k_2N(1+N^\frac{1}{2})\sum_{i=1}^Nm_i\mathrm{e}^{-a_it}.
\end{align*}

Considering the event-triggered scheme \eqref{trigger_scheme}, when $t<t_{k+1}^i$, we have $G_i(t)\leq 0$, which results in $\sum_{i=1}^N\|\varepsilon_i\|\leq \sum_{i=1}^Nm_i\mathrm{e}^{-a_it}$ and thus
\begin{equation}\label{error_inequality}
\begin{aligned}
    \sum_{i=1}^N\sum_{j\in\mathcal{N}_i}\|\varepsilon_i-\varepsilon_j\|\leq&
\sum_{i=1}^N\sum_{j=1}^N\|\varepsilon_i-\varepsilon_j\|\\
\leq& \sum_{i=1}^N\sum_{j=1}^N\|\varepsilon_i\|+\sum_{i=1}^N\sum_{j=1}^N\|\varepsilon_j\|\\
    =&2N\sum_{i=1}^N\|\varepsilon_i\|
    \leq
    2N\sum_{i=1}^Nm_i\mathrm{e}^{-a_it}.
\end{aligned}
\end{equation}

Therefore, we can obtain that
\begin{align*}
\dot{V}_{2}\leq&-\frac{k_2}{2}\sum_{i=1}^N\sum_{j\in\mathcal{N}_i}\|x_i-x_j\|
\leq0.
\end{align*}

According to the Barbalat Lemma \cite{slotine1991applied}, $\widehat{x}_i$ will converge to 0 as $t\rightarrow\infty$. In other words, $\|x_i-x_j\|\rightarrow 0, \forall i,j\in \mathcal{V}$ as $t\rightarrow\infty$.
Based on Lemma \ref{le_op_fix}, we can deduce that all the state $x_i$ will eventually converge to the optimal trajectory $x^*(t)$ of the DTOP \eqref{DTOP}.
$\hfill{} \Box$
\end{pf}

% Note that the Zeno behavior should be avoided in event-triggering approaches. In the following result confirms that the Zeno behavior can be avoided in our scheme.
\begin{lemma}\label{le_zeno}
Under the conditions in Theorem \ref{th_fixtop},
     the Zeno behavior is avoided by the proposed scheme with the triggering function \eqref{trigger_function}.
\end{lemma}
\begin{pf}
     The Lemma is proved by contradiction. Assume that the Zeno behavior occurs, then there exists a constant $T^*>0$ such that $\lim_{k\rightarrow \infty}t_k^i=T^*$.

     For any $t\in[t_k^i,t_{k+1}^i)$, the time derivative of $\varepsilon_i$ defined in \eqref{measurement_error} can be obtained as:
\begin{align*}
\dot{\varepsilon}_i=&-\dot{x}_i\\=
&k_1s_i+k_2\sum_{j\in \mathcal{N}_i}\mathrm{sign}(\widetilde{x}_i-\widetilde{x}_j)+k_3\mathrm{sign}(\nabla^2_{x_i}f(x_i,t)s_i).
\end{align*}
It follows from Lemma \ref{le_op_fix} that
$s_i$, for all $ i\in\mathcal{V}$, is bounded, and moreover, $\| \mathrm{sign} (\cdot) \|$ is also bounded for any vector. As a result,  $\| \dot{\varepsilon}_i \|$ is bounded, and we can denote $\|\dot{\varepsilon}_i\|\leq Q$ with $Q>0$. Then we have
$$ \|\varepsilon_i (t)\|= \left\|\int_{t_k^i}^t\dot{\varepsilon}_i\mathrm{d}\tau
\right\| \leq\int_{t_k^i}^t\|\dot{\varepsilon}_i\|\mathrm{d}\tau\leq Q(t-t_k^i),$$
for any $t\in[t_k^i,t_{k+1}^i)$.
% Therefore,
% \begin{align*}
% ||\dot{\varepsilon}_i||\leq&k_2||\varepsilon_i||
% +k_1||s_i||\\
% &+||k_3\mathrm{sign}(\nabla^2_{x_i}f(x_i,t)s_i)||
% +||k_4\sum_{j\in \mathcal{N}_i}\mathrm{sign}(\widetilde{x}_i-\widetilde{x}_j)||\\
% \leq&k_2||\varepsilon_i||+Q.
% \end{align*}
% Thus, we have
% \begin{align*}
%     ||\varepsilon_i||\leq\frac{Q}{k_2}\big(\mathrm{e}^{k_2(t-t_k^i)}-1\big),~~\forall t\in [t_k^i,t_{k+1}^i).
% \end{align*}
Then, considering the event-triggering scheme \eqref{trigger_scheme}, we have
\begin{align*}
Q(t_{k+1}^i-t_k^i)&\geq\|\varepsilon_i(t_{k+1}^i)\|
>m_i\mathrm{e}^{-a_it_{k+1}^i}.
\end{align*}
Note that $\lim_{k\rightarrow \infty}t_k^i=T^*$ implies $\lim_{k\rightarrow \infty} m_i\mathrm{e}^{-a_it_{k+1}^i} = m_i\mathrm{e}^{-a_i T^*}>0$, while $\lim_{k\rightarrow \infty}(t_{k+1}^i - t_k^i)=0$. This leads to a contradiction.

Therefore, there cannot exist a constant $T^*>0$ such that $\lim_{k\rightarrow \infty}t_k^i=T^*$. As a result, the Zeno behavior does not occur.
$\hfill{} \Box$
\end{pf}

\begin{remark}
In Theorem \ref{th_fixtop} and Lemma \ref{le_op_fix}, the assumption that $\nabla^2_{x_i}f_i(x_i,t)=\nabla^2_{x_j}f_j(x_j,t)$ holds for all $i,j\in \mathcal{V}, t\geq 0$  is strict, but it is common in distributed time-varying optimization problems  \cite{10402058,8715378,8100702,7518617}. As our future research, we will explore  possible approach
in  distributed time-varying optimization to remove this assumption.
\end{remark}

\section{Numerical Simulation}\label{sec4}
In this section,   a case study about a state-of-charge (SOC) balancing problem for a battery energy storage system (BESS) is presented to show the effectiveness of the proposed method of Section \ref{sec3}.

%\begin{example}\label{soc}
%We apply the proposed  approach to a package-level state-of-charge (SOC) balancing problem  for a battery energy storage system (BESS).
     BESSs are becoming more important in managing power systems, for instance, to improve the power quality of renewable-energy hybrid power generation systems
 %     and to satisfy complex
 % energy requirements  with flexible energy management solutions
 \cite{6473871}. By releasing power at peak times and storing power at off-peak times, a BESS system can help to reduce the increasing electric demand and drive power generation to run at optimal efficiency \cite{7548310}.
A general BESS structure consists of multiple battery packages, and each
battery package owes a group of  battery cells. Each battery package is equipped with a battery management system which can monitor and balance the SOCs of all the cells~\cite{7548310}. The package-level SOC balancing can protect battery packages from overcharging or discharging and stabilize grid frequency and voltage \cite{7548310}.

Consider a BESS with $N$ battery packages with  $P_i\in\mathbb{R}$  the power output of the $i$th battery package, where $P_i>0$ during discharging operations and $P_i<0$ during charging operations.
% Generally, bidirectional ac-dc converters are applied to connect the batteries with the grid \cite{7548310}. Let $Q$  be the total capacity of each battery package and  $V_d$ a battery related constant. When the grid-connected battery package system is in steady state,  in the dc side, the output current of the $i$th battery $I_{bi}\approx I_{oi}$ with $I_{oi}$ being the input current of the converter and the voltage of the  $i$th converter can be assumed to be the same for the same type of battery packages. Then, for the $i$th battery package,  the connection between $P_i$ and $SOC_i$ is
% $$\dot{SOC}_i \approx -\frac{E_i}{QV_d}P_i,
% $$ where $E_i$ is influenced by the Coulombic efficiency \cite{7548310,8100702}.
Let $P^*(t)$ be the desired power output for the whole BESS,
 % where $P^*(t)>0$ for the discharging operation and $P^*(t)<0$ for the charging operation.
then in the package-level SOC balancing problem, we need to regulate each $P_i$ such that the total power output $P(t): =\sum_{i=1}^NP_i(t)$ satisfies
\begin{align}\label{soc_goal}
    \limsup_{t\rightarrow \infty}|P(t)-P^*(t)|<\delta,
\end{align}
for some small $\delta>0$, and further to ensure that the
state-of-charge of all the battery packages are nearly equal, see more details in \cite{7548310}.
% $$\limsup_{t\rightarrow \infty}|SOC_i(t)-SOC_j(t)|<\delta$$ for some small $\delta>0$
The BESS can be regarded as a multi-agent system, and a distributed time-varying optimization method can be applied to achieve \eqref{soc_goal}, where every package can track its own desired power output according to its local optimal strategy.

% Similar to \cite{8100702}, we add a distributed optimization process to obtain $P^*(t)$ rather than a global command generator as in \cite{7548310}. In this way, every package can set its desired power output according to its own local optimal strategy.
% A block diagram of the package-level balancing scheme is shown in Fig. \ref{fig_scenario}, where $I_i$ is the current of the $i$th converter, $v$ is the voltage of the distribution bus and BMS denotes the battery management system.
% \begin{figure}
%     \centering
%        % \includegraphics[height=4.5cm]{figures/scenario.pdf}
% \includegraphics[height=5cm]{Automatica/figures/block_soc.pdf}
%         \caption{Block diagram of the package-level balancing scheme. It is inspired by Fig. 2 in \cite{7548310} and Fig.1 in \cite{8100702}.}
%         \label{fig_scenario}
% \end{figure}

% In what follows, we mainly aim at obtaining  $P^*(t)$ by the distributed optimization approach \eqref{controller}.
\begin{figure}
    \centering
       \includegraphics[height=6cm]{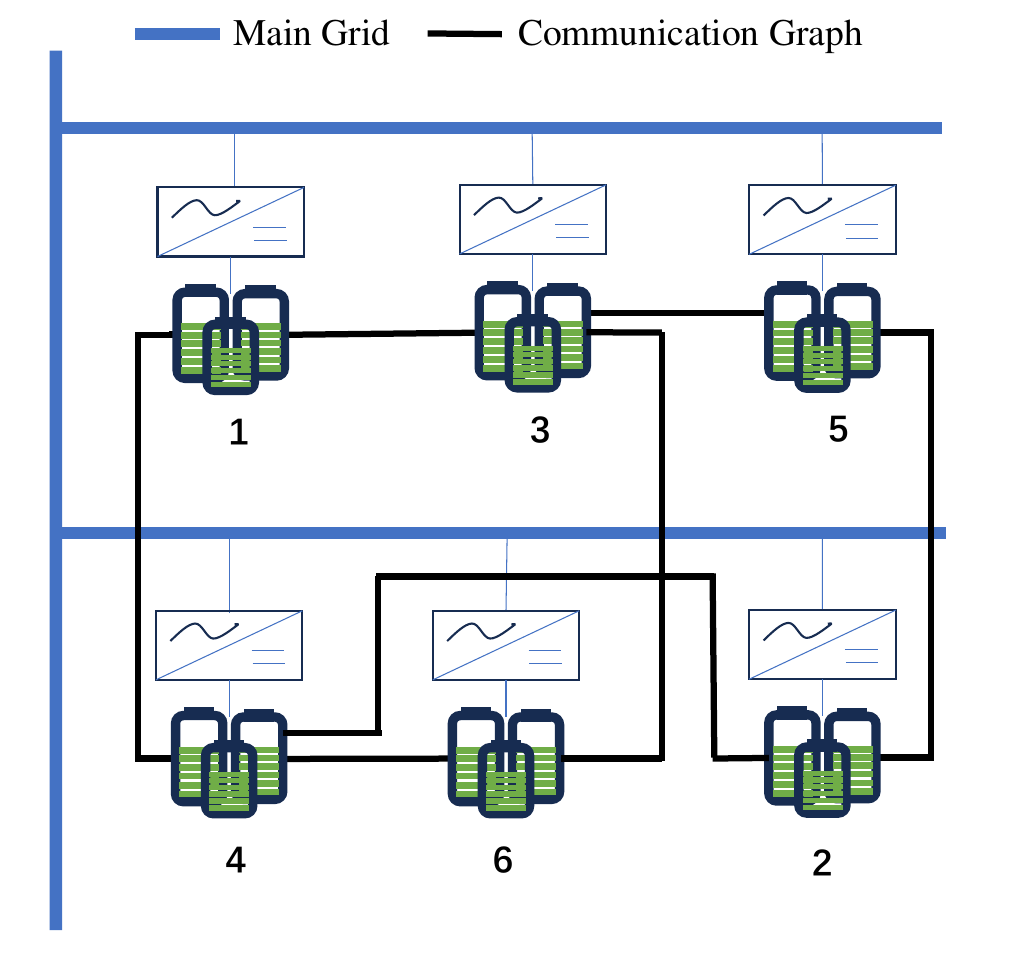}
        \caption{The communication graph among battery packages. }
        \label{fig_toplogy}
\end{figure}
Note that the goal \eqref{soc_goal} is achieved when optimizing $P_i(t)$ to satisfy $P_i(t)=P_{\mathrm{ave}}^*(t)$, where
$$ P_{\mathrm{ave}}^*(t) : = \frac{1}{N}P^*(t).$$
Consider six battery packages ($N=6$), and the communication graph among battery packages is given in Fig.~\ref{fig_toplogy}.
Then, we formulate the charging problem in the form of DTOP \eqref{DTOP} as
% Denote $P_{avg}^*(t)=\frac{1}{N}P^*(t)$. If $P_i$ is controlled such that $P_i=P^*_{avg}$, the goal \eqref{soc_goal} is achieved.
% % Without loss of generality, denote $P^*(t)=\sum_{i=1}^NP^*_i(t)$.
% From a distributed optimization perspective, the goal of obtaining $P^*(t)$
% can be abstracted
%  The problem of obtaining $P^*(t)$ is rewritten as:
\begin{equation} \label{soc}
\begin{aligned}
    \min_{P_i(t)} ~~&\sum_{i=1}^6  f_i(P_i(t), t)\\
\mathrm{s.t.}~~~&P_i(t)=P_j(t),~~ \forall i,j=1,2,\dots,N,
\end{aligned}
\end{equation}
%  $P^*_1(t)=t,P^*_2(t)=0.2t,P^*_3(t)=0.5\sin(2t),0.5\sin(2t),P^*_4(t)=0.1\cos(2t),P^*_5(t)=0.5t,P^*_6(t)=1.2t$
where the local objective functions for each package as
\begin{equation}\label{example_soc}
\begin{split}
    &f_1(P_1(t),t)=\frac{1}{2}(P_1(t)-t)^2,\\
     &f_2(P_2(t),t)=\frac{1}{2}(P_2(t)-0.2t)^2,\\
      &f_3(P_3(t),t)=\frac{1}{2}(P_3(t)-0.5\sin(2t))^2,\\
      &f_4(P_4(t),t)=\frac{1}{2}(P_4(t)-0.1\cos(2t))^2,\\
      &f_5(P_5(t),t)=\frac{1}{2}(P_5(t)-0.5t)^2,\\
      &f_6(P_6(t),t)=\frac{1}{2}(P_6(t)-1.2t)^2.
      \end{split}
\end{equation}
It can be verified that Assumptions \ref{as3} and \ref{as2} are satisfied for each $f_i(P_i(t), t)$ with $\overline{f} = 1.2$ and $\underline{\lambda}_i=1$. Also, the condition on the Hessian matrix for each objection function fulfills \eqref{hessian}.

We implement the distributed optimization approach~\eqref{controller} to steer all $P_i(t)$ to track the average of the desired power output $P_{\mathrm{ave}}^*(t)$ in a distributed way. Here, intermittent communication is achieved by the event-triggered distributed controller \eqref{controller}. In the simulation, we set $k_1=5, k_2=1, k_3=15$, where $k_3$ satisfies the condition \eqref{K3}.  Moreover, the tuning parameters in the triggering functions of the agent are set as $a_1=0.9,a_2=0.7,a_3=0.9,a_4=0.9,a_5=0.7,a_6=0.7$, $m_1=3,m_2=2,m_3=3,m_4=3,m_5=2,m_6=2$.

\begin{figure} [!h]
\centering
\includegraphics[height=6cm]{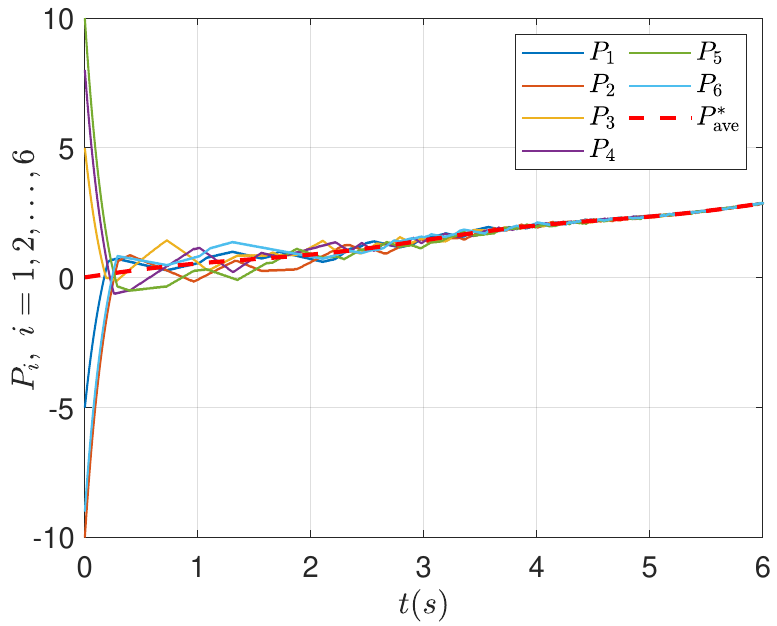}
        \caption{Trajectories of power output. }
        \label{fig_state_soc}
\end{figure}

\begin{figure} [!h]
     \centering
        \includegraphics[height=6cm]{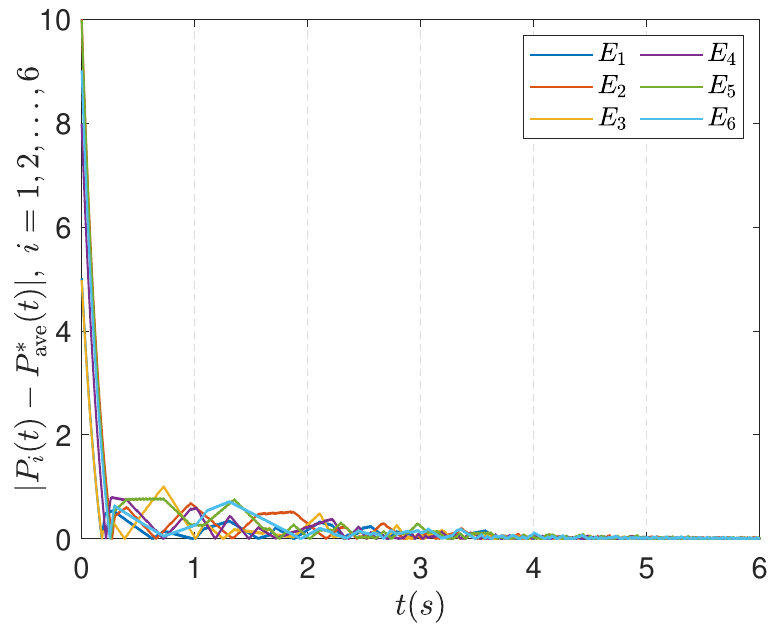}
        \caption{The errors between $P_i$ and $P^*$. }
        \label{fig_errors_soc}
\end{figure}

\begin{figure}
\centering
        \includegraphics[height=6cm]{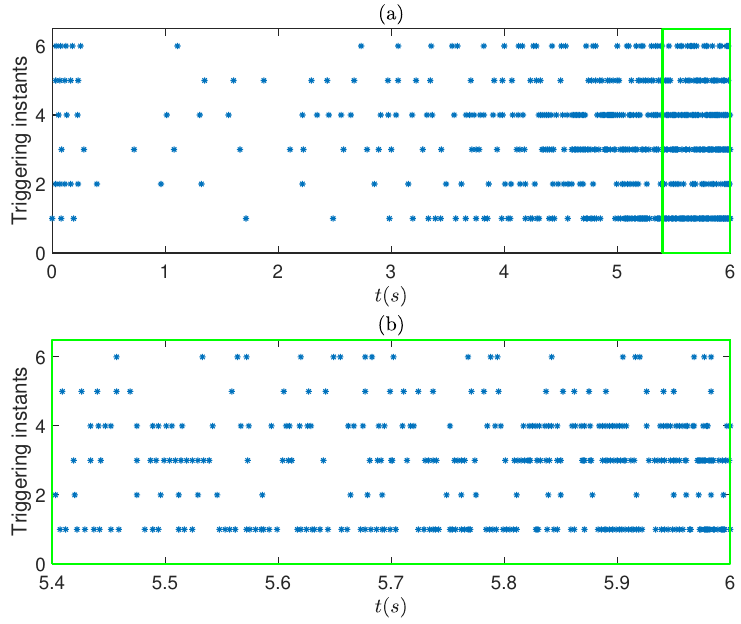}
        \caption{The triggering instants of each battery package. (a) The communication instants triggered in the full interval $t\in[0, 6]$. (b) Zoomed-in figure showing the communication instants in the region $t\in[5.4, 6]$ (boxed area in (a)).}
        \label{fig_triggering_instants}
\end{figure}

 The trajectories of $P_i(t)$ generated by the distributed controller \eqref{controller} are shown in Fig.~\ref{fig_state_soc}. We observe that these trajectories can converge to the optimal reference trajectory $P_{\mathrm{ave}}^*(t)$.
 % Here, we obtain the optimal trajectory the by optimal condition (Lemma 1 in \cite{9855238}).
 Furthermore, the error between the power trajectories and the desired power output is illustrated in Fig.~\ref{fig_errors_soc}, which shows a clear delay over a short time.
 % We observe that the errors fluctuate around zero with small amplitudes because the communication between battery packages is intermittent and the tracking error may gradually increase when no communication.
In addition, the triggering instants of each agent are recorded in Fig.~\ref{fig_triggering_instants}, which shows that only intermittent communication is required for each agent following the event trigger scheme in \eqref{trigger_scheme}. Therefore, compared to the algorithms in  \cite{10339857,10402058,8100702,9855238}, where
 continuous-time communication is required, the proposed method with the distributed controller \eqref{controller} can significantly reduce the communication cost in the distributed implementation.


\begin{figure}[!h]
\centering
        \includegraphics[height=6cm]{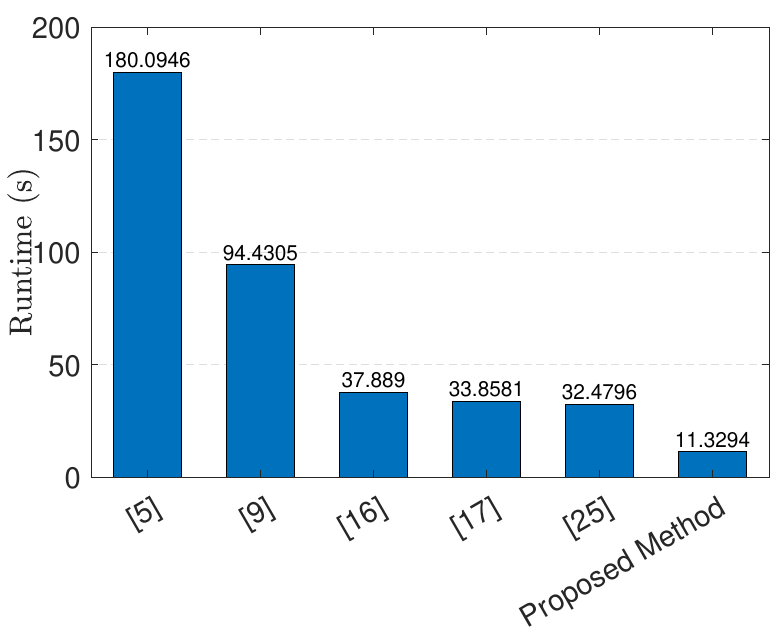}
        \caption{Comparison of the runtime for solving the problem \eqref{soc} using the different methods for \cite{10339857,10402058,doi:10.1080/00207721.2020.1801885,8100702,9855238}. }
        \label{fig_runtime}
\end{figure}

% \begin{figure}[!h]
% \centering
%         % \includegraphics[height=6cm]{figures/triggering_instants.pdf}
%         \includegraphics[height=6cm]{Automatica/figures/detail_instants_soc.pdf}
%         \caption{The partial triggering instants of each battery package. }
%         \label{fig_detail_instants}
% \end{figure}
Next, we compare the runtime of the distributed controller \eqref{controller} and the algorithms in~\cite{10339857,10402058,doi:10.1080/00207721.2020.1801885,8100702,9855238} when solving the problem \eqref{soc} over the time interval $t \in [0,6]$. The results are presented in Fig.~\ref{fig_runtime}, which shows that the runtime of our method is notably lower than that of the referenced algorithms. This improvement is because the distributed controller \eqref{controller} does not need to compute the inverse of the Hessian matrix of the local objective function in real time, and it also does not need to continuously exchange information on the partial derivatives of the local objective function in real time among the agents. As a result, compared to the algorithms in \cite{10339857,10402058,doi:10.1080/00207721.2020.1801885,8100702,9855238}, the proposed method can reduce the computational cost for solving the DTOP.

% \begin{table*}[ht]
% \centering
% \caption{Runtime for solving the problem \eqref{example_soc} with different methods.}
% \label{tab_runtime}
% \begin{tabular}{ccccccc}
% \toprule[1pt]
% Methods & \cite{9855238} & \cite{8100702} & \cite{10402058} & \cite{10339857} & \cite{doi:10.1080/00207721.2020.1801885} & Proposed Method\\
% \midrule[1pt]
% Runtime (s)& 32.479574 & 33.858098 & 94.430525 & 180.094569 & 37.888969 & 11.329429\\
% \bottomrule[1pt]
% \end{tabular}
% \end{table*}
%\end{example}

\section{Conclusion}\label{sec5}
In this paper, we have proposed a new method for solving a class of distributed time-varying optimization problems by integrating an event-triggered scheme into a distributed neurodynamic optimization approach. The event-triggered scheme has been shown to be effective in reducing communication frequency, thereby conserving communication resources.
Furthermore, the proposed distributed controller does not need to compute the inverse of the Hessian of the local objective functions, contributing to a low computational cost. A case study on a battery charging problem has shown the effectiveness of the distributed neurodynamic optimization approach.
% in complex computational problem and engineering application.
In terms of future research, we plan to further generalize the method, which does not require identical Hessian matrices of the local objective functions to solve distributed time-varying optimization problems.

\bibliographystyle{abbrv}%automatic
\bibliography{event_triggered,I}

\end{document}